\documentclass[preprint]{elsarticle}

\usepackage{amssymb}
\usepackage{amsfonts}
\usepackage{amsmath}
\usepackage{graphicx}
\usepackage{amsthm}
\usepackage{algorithm}
\usepackage{epic}

\newtheorem{rem}{Remark}

\newtheorem{exmpl}{Example}
\newenvironment{example}{\begin{exmpl}\rm}{\end{exmpl}}

\def\diag{\mathop{\mathrm{diag}}}

\begin{document}

\begin{frontmatter}

\title{Forward stable eigenvalue decomposition of rank-one
  modifications of diagonal matrices}
\author[FESB]{N.~Jakov\v{c}evi\'{c}~Stor\corref{c1}\fnref{f1}}
\ead{nevena@fesb.hr}
\author[FESB]{I.~Slapni\v{c}ar\fnref{f1}}
\ead{ivan.slapnicar@fesb.hr}
\author[PEN]{J.~L.~Barlow\fnref{f2}}
\ead{barlow@cse.psu.edu}

\cortext[c1]{Corresponding author}
\fntext[f1]{The research of Ivan Slapni\v{c}ar and Nevena Jakov\v{c}evi\'{c}
  Stor was supported by the Ministry of Science, Education and
  Sports of the Republic of Croatia under grant 023-0372783-1289.}
\fntext[f2]{The research of Jesse L. Barlow  was supported
by the National Science Foundation under grant  CCF-1115704.}

\address[FESB]{Faculty of Electrical Engineering, Mechanical Engineering and
Naval Architecture, University of Split, Rudjera Bo\v{s}kovi\'{c}a
32, 21000 Split, Croatia}

\address[PEN]{Department of Computer Science and Engineering, The Pennsylvania State
University, University Park, PA 16802-6822, USA}

\begin{abstract}
We present a new algorithm for solving an eigenvalue problem for a real
symmetric matrix which is a rank-one modification of a diagonal matrix.
The algorithm computes
each eigenvalue and all components of the corresponding eigenvector with
high relative accuracy in $O(n)$ operations.
The algorithm is based on a shift-and-invert approach.
Only a single element of the inverse of the shifted matrix eventually needs to
be computed with double the working precision.
 Each eigenvalue and the corresponding
eigenvector can be
computed separately, which makes the algorithm adaptable for parallel
computing. Our results extend to the complex Hermitian case. The algorithm is
similar to the algorithm for solving the eigenvalue problem for real symmetric arrowhead
matrices from: N.~Jakov\v{c}evi\'{c}~Stor, I.~Slapni\v{c}ar and J.~L.~Barlow,
{Accurate eigenvalue decomposition of real symmetric arrowhead matrices and
  applications}, Lin.\ Alg.\ Appl., 464 (2015).

\end{abstract}


\end{frontmatter}

\section{Introduction and Preliminaries}

In this paper we consider the eigenvalue problem for an $n\times n$ real
symmetric matrix $A$ of the form

\begin{equation}
A= D +\rho z z^T,
\label{A}
\end{equation}%
where
\begin{equation*}
D=\mathop{\mathrm{diag}}(d_{1},d_{2},\ldots ,d_{n})\text{ }
\end{equation*}%
is a diagonal matrix of order $n$,%
\begin{equation*}
z=\left[
\begin{array}{cccc}
\zeta _{1} & \zeta _{2} & \cdots & \zeta _{n}%
\end{array}%
\right] ^{T}  \label{z}
\end{equation*}%
is a vector and $\rho \neq 0$ is a scalar.
Notice that $A$ is a rank-one modification of a diagonal matrix. Subsequently,
we shall refer to such matrices as ``diagonal-plus-rank-one'' (DPR1) matrices.
DPR1 matrices arise, for example, in solving symmetric real tridiagonal
eigenvalue problems with the divide-and-conquer method \cite{Cup81},
\cite{DS87}, \cite{GE95}, 
\cite[Sections 3.2.1 and 3.2.2]{Ste01}, \cite[Section III.10]{VVM08}.

Without loss of generality, we make the following assumptions:
\begin{itemize}
\item[-] $\rho>0$ (otherwise we consider the matrix $A=-D-\rho zz^T$),
\item[-] $A$ is irreducible, that is,
$ \zeta _{i}\neq 0, i=1,\ldots, n$, and
$
d_{i}\neq d_{j},\text{ for all }i\neq j,\text{ }i,j=1,\ldots ,n
$, and
\item[-]
the diagonal elements of $D$ are decreasingly ordered,
\begin{equation}
d_{1}>d_{2}>\cdots >d_{n}.  \label{order}
\end{equation}%
\end{itemize}
Indeed, if $\zeta _{i}=0$ for some $i$, then the diagonal element $d_{i}$ is an
eigenvalue whose corresponding eigenvector is the $i$-th unit vector, and
if $d_{i}=d_{j}$, then $d_{i}$ is an
eigenvalue of the matrix $A$ (we can reduce the size of the problem by
annihilating $\zeta_j$ with a Givens rotation in the $(i,j)$-plane).
Ordering of the diagonal elements of $D$ is attained by symmetric row and
column pivoting.

Let
\begin{equation*}
A=V\Lambda V^{T}  \label{Aeigendec}
\end{equation*}%
be the eigenvalue decomposition of $A$, where
\begin{equation*}
\Lambda =\mathop{\mathrm{diag}}(\lambda _{1},\lambda _{2},\ldots ,\lambda
_{n})
\end{equation*}%
is a diagonal matrix whose diagonal elements are the eigenvalues of $A$, and
\begin{equation*}
V=\left[
\begin{array}{ccc}
v_{1} & \cdots & v_{n}%
\end{array}%
\right]
\end{equation*}%
is an orthonormal matrix whose columns are the corresponding eigenvectors.

The eigenvalue problem for a DPR1 matrix $A$ can be solved by any of the
standard methods for the symmetric eigenvalue problem (see, for example
\cite{Wil65,Par80}).
However, because of the special structure of diagonal-plus-rank-one matrices, we
can use the following approach.
The eigenvalues of $A$ are the zeros of the secular function (see, for
example, \cite{Cup81}  and \cite[Section 8.5.3]{GV96}):
\begin{equation}
f(\lambda )=1+\rho\sum_{i=1}^{n}\frac{\zeta _{i}^{2}}{%
d_{i}-\lambda }=1 +\rho z^{T}(D-\lambda I)^{-1}z,  \label{Pick}
\end{equation}%
and the corresponding eigenvectors are given by
\begin{equation}
v_{i}=\frac{x_{i}}{\left\Vert x_{i}\right\Vert _{2}},\text{ \ \ \ }
x_{i}=( D-\lambda _{i}I) ^{-1}z,\quad i=1,\ldots ,n.  \label{Aeigenvec}
\end{equation}%
Diagonal elements of the matrix $D$, $d_{i}$, are called poles of the
function $f$. It is easy to see that, for $\rho>0$,
$f$ is strictly increasing between the poles,
implying the strict interlacing property
\begin{equation}
\lambda _{1}>d_{1}>\lambda _{2}>d_{2}>\cdots >
\lambda _{n}> d_n.  \label{interlace}
\end{equation}

The formulae (\ref{Pick}) and (\ref{Aeigenvec}) are simple,
and have been used to solve similar eigenvalue problems
\cite{Bar93, BN78, Cup81, DS87}.
but maintaining orthogonality among the eigenvectors $v_i$
requires all of
the eigenvalues $\lambda_i$ to be computed with high accuracy \cite{GE95}.
In other words, if the computed eigenvalues
are not accurate enough,
then the computed eigenvectors
may not be sufficiently orthogonal
(see Example \ref{example3}). The existing algorithms for DPR1 matrices
\cite{Cup81,DS87,GE95} obtain orthogonal eigenvectors with the following
procedure:

\begin{itemize}
\item[-] compute the eigenvalues $\tilde \lambda_i$ of $A$ by solving (\ref%
{Pick}),

\item[-] construct a new matrix
\begin{equation*}
\label{eq:Gumatrix}
\tilde{A}=
D + \rho \tilde{z}
\tilde{z}^{T}
\end{equation*}%
by solving an inverse problem with the prescribed eigenvalues,
\item[-] compute the eigenvectors of $\tilde{A}$ by (\ref{Aeigenvec}) but
  using $\tilde z$ instead of $z$.
\end{itemize}
The eigenvectors computed by this algorithm are
orthogonal to machine precision (for details see
\cite{GE95,Cup81,DS87,Bar93}).
This results in an algorithm which requires only $O(n^{2})$ computations and
$O(n)$ storage for eigenvalues and $O(n)$ storage for each eigenvector.
This algorithm is implemented in the LAPACK
subroutine \verb!DLAED9! and its subroutines \cite{ABB99}.

Our algorithm uses a different approach and is forward stable, that is, it
computes all eigenvalues and all individual components of the
corresponding eigenvectors of a given arrowhead matrix of floating-point
numbers to almost full accuracy, a feature which no other method has.
The accuracy of the eigenvectors and
their numerical orthogonality follows from the high relative accuracy of the
computed eigenvalues. Each eigenvalue and the corresponding eigenvector is
computed independently of the others in $O(n)$ operations, making our
algorithm suitable for parallel computing.

The algorithm is based on a shift-and-invert technique. Basically, an
eigenvalue $\lambda $ is computed from the largest or the smallest eigenvalue
of the inverse of the matrix shifted to the pole $d_{i}$ which is nearest to
$\lambda $, that is,
\begin{equation}  \label{8}
\lambda=\frac{1}{\nu}+d_i,
\end{equation}
where $\nu$ is either largest or smallest eigenvalue of the
matrix
\begin{equation*}
A_i^{-1}\equiv (A-d_iI)^{-1}.
\end{equation*}

The algorithm and its error analysis are similar to the algorithm for
arrowhead matrices from  \cite{JSB13}, thus, the present paper can be viewed as a
note related to \cite{JSB13}.

The organization of the paper is the following.
In Section 2, we describe our algorithm named $dpr1eig$ and give error bounds.
We also discuss fast secular equation solvers and three implementations of
the double the working precision. In Section 3, we illustrate our algorithm
with few examples.

\section{The algorithm}

Let $A$ be an irreducible DPR1 matrix of the form (\ref{A}), with the diagonal
elements of $D$ ordered as in (\ref{order}), and $\rho>0$.
Let $\lambda $ be an eigenvalue of $A$, let $v$ be its eigenvector, and let $%
x$ be the unnormalized version of $v$ from (\ref{Aeigenvec}). Let $d_{i}$ be
a pole which is closest to $\lambda $. Clearly, from (\ref{interlace}) it
follows that either $\lambda =\lambda _{i} $ or $\lambda =\lambda _{i+1}$.
Let $A_{i}$ be the shifted matrix
\begin{align*}
A_{i}=A-d_{i}I=\left[
\begin{array}{ccc}
D_{1} & 0 & 0 \\
0 & 0 & 0  \\
0 & 0 & D_{2}
\end{array}%
\right]
+\rho \begin{bmatrix} z_{1} \\ \zeta _{i} \\ z_{2}
\end{bmatrix}
\begin{bmatrix}
z_{1}^{T} & \zeta _{i} & z_{2}^{T}
\end{bmatrix},
\label{Ai}
\end{align*}%
where%
\begin{align*}
D_{1}& =\mathop{\mathrm{diag}}(d_{1}-d_{i},\ldots ,d_{i-1}-d_{i}), \\
D_{2}& =\mathop{\mathrm{diag}}(d_{i+1}-d_{i},\ldots ,d_{n}-d_{i}), \\
z_{1}& =\left[
\begin{array}{cccc}
\zeta _{1} & \zeta _{2} & \cdots & \zeta _{i-1}%
\end{array}%
\right] ^{T}, \\
z_{2}& =\left[
\begin{array}{cccc}
\zeta _{i+1} & \zeta _{i+2} & \cdots & \zeta _{n}%
\end{array}%
\right] ^{T}.
\end{align*}%
Notice that $D_{1}$ $(D_{2})$ is positive (negative) definite.

Obviously, $\lambda $ is an eigenvalue of $A$ if and only if
\begin{equation*}
\mu =\lambda -d_{i}
\end{equation*}%
is an eigenvalue of $A_{i}$, and they share the same
eigenvector. 

The inverse of $A_{i}$ is a permuted arrowhead matrix
\begin{equation}
A_{i}^{-1}=\left[
\begin{array}{ccc}
D_{1}^{-1} & w_{1} & 0 \\
w_{1}^{T} & b & w_{2}^{T} \\
0 & w_{2} & D_{2}^{-1}
\end{array}%
\right] ,  \label{invAi}
\end{equation}%
where%
\begin{align}
w_{1}& =-D_{1}^{-1}z_{1}\frac{1}{\zeta _{i}},  \notag \\
w_{2}& =-D_{2}^{-1}z_{2}\frac{1}{\zeta _{i}},  \notag \\
b& =\frac{1}{\zeta _{i}^{2}}\left(
\frac{1}{\rho}+z_{1}^{T}D_{1}^{-1}z_{1}+z_{2}^{T}D_{2}^{-1}z_{2}\right).
\label{b_ob}
\end{align}
The above formulas for the inverse, which can be verified directly, can also be deduced from \cite[Fact
2.16.4]{Ber09}, \cite[pp.\ 225]{DeVB06} or \cite[Theorem 1]{ElRo81}.  
The computation of the scalar $b$ in (\ref{b_ob}),
is critical to how well we are able to compute $\lambda$.

The eigenvalue $\nu$ of a real symmetric arrowhead matrix $A_i^{-1}$
from (\ref{invAi})
is a zero of the secular equation
(see, for example \cite{OlSt90,JSB13})
\begin{equation}\label{secAi}
g(\nu)=b-\nu-w^T (\Delta -\nu I)^{-1}w=0,
\end{equation}
where
$$
\Delta=\begin{bmatrix}D_1 & \\ & D_2
\end{bmatrix}, \qquad
w=\begin{bmatrix}w_1 \\ w_2
\end{bmatrix}.
$$
Once $\nu$ is computed, we compute $\mu=1/\nu$.
The normalized and unnormalized eigenvectors $v$ and $x$ are computed by
applying (\ref{Aeigenvec}) to the matrix $A_i$, that is,
\begin{equation}
x=\begin{bmatrix} x_1\\ \vdots \\ x_n
\end{bmatrix}
=
\left[
\begin{array}{c}
\left( D_{1}-\mu I\right) ^{-1}z_{1} \\
-\displaystyle\frac{\zeta _{i}}{\mu } \\
\left( D_{2}-\mu I\right) ^{-1}z_{2}
\end{array}
\right], \quad v=\frac{x}{\| x\|_2}.  \label{eigenvecAi}
\end{equation}

If $\lambda $ is an eigenvalue of $A$ which is closest to the pole $d_{i}$,
then $\mu $ is the eigenvalue of matrix $A_{i}$ which is closest to zero and
\begin{equation*}
\nu=\frac{1}{\mu}=\pm \left\Vert A_{i}^{-1}\right\Vert _{2}.
\end{equation*}
We say that $\nu$ is the {\em largest absolute eigenvalue} of $A_{i}^{-1}$.
In this case, if all entries of $A_{i}^{-1}$ are computed with high relative
accuracy, then, according to standard perturbation theory, any reasonable
algorithm can compute $\nu$
to high relative accuracy (see Section \ref{sec:fast}).

Throughout the paper, we assume that the computations are carried out in the 
standard floating-point arithmetic with the machine
precision $\varepsilon_M= 2^{-52}\approx 2.2204\cdot 10^{-16}$ (see
\cite[Chapter 2]{Hig96} for details). Thus, the floating-point numbers
have approximately 16 significant decimal digits. The term ``double the
working precision'' means that the computations are performed with numbers
having approximately 32 significant decimal digits, or with the machine
precision $\varepsilon_M^2$ or smaller.

Notice that all entries of $A_{i}^{-1}$ are computed to high relative
accuracy using standard precision, except possibly $b$ in (\ref{b_ob}).
For example, using the standard model from \cite[Section 2.2]{Hig96}, the
error analysis for the respective indices $k$ gives
\begin{align*}
fl([D_1]_k)&=\frac{1}{d_k-d_1}(1+\varepsilon_1), \quad |\varepsilon_1|\leq 2\varepsilon_M,\\
fl([w_1]_k)&=\frac{\zeta_k}{\zeta_i(d_k-d_i)}(1+\varepsilon_2), \quad |\varepsilon_2|\leq 3\varepsilon_M,
\end{align*}

If $b$ is not computed to
high relative accuracy and it influences $\left\Vert A_{i}^{-1}\right\Vert
_{2}$, it is sufficient to compute it with double the working precision.
Whether double the working precision is needed is determined as follows:
set
\begin{align}
K_{b}& =\frac{1 +\rho
z_{1}^{T}D_{1}^{-1}z_{1} - \rho
z_{2}^{T}D_{2}^{-1}z_{2} }{\left\vert
1+\rho z_{1}^{T}D_{1}^{-1}z_{1}+\rho z_{2}^{T}D_{2}^{-1}z_{2}\right\vert },\nonumber
\\
K_z&= \frac{1}{%
|\zeta _{i}|}\sum\limits_{\substack{ j=1  \\ j\neq i}}^{n}|\zeta
_{j}|, \nonumber\\
\kappa _{\nu }& \leq
\min\big\{(n+4)\sqrt{n}\, K_{b},
 3\sqrt{n}+(n+4)\big(1+ 2 K_z \big)\big\}.
\label{knu2}
\end{align}
Here $K_b$ measures whether $b$ is computed with high relative accuracy,
$K_z$ measures whether $b$ influences $\left\Vert A_{i}^{-1}\right\Vert
_{2}$, and $\kappa_\nu$ measures the accuracy of the exact eigenvalue
$\widehat \nu$ of
the computed matrix $fl(A_i^{-1})$,
$$
\widehat \nu =\nu(1+\kappa_\nu \varepsilon_M),
$$
similarly as in \cite[Theorem 5]{JSB13}.

If $\kappa_\nu\gg O(n)$, then $b$ needs to be computed in double the
working precision (see section \ref{double}).
The details of the proofs of the above facts are similar to the proofs of
\cite[Theorems 5 and 7]{JSB13}.

If $\lambda $ is an eigenvalue of $A$ which is not closest to the pole
$d_{i}$, then $\mu $ is not the eigenvalue of $A_{i}$ which is closest to
zero. Further, $\left\vert \nu \right\vert <\left\Vert A_{i}^{-1}\right\Vert
_{2} $, and the quantity
\begin{equation}
K_{\nu} = \frac{\left\Vert A_{i}^{-1}\right\Vert _{2}}{\left\vert \nu
\right\vert}  \label{K0}
\end{equation}%
tells us how far $\nu $ is from the largest absolute eigenvalue of $%
A_{i}^{-1}$. If $K_{\nu}\gg 1$, then the standard perturbation theory does
not guarantee that the eigenvalue $\mu $ will be computed with high relative
accuracy. One remedy to this situation is to use non-standard shifting as
follows:
\begin{enumerate}
\item[(R1)]
we can compute $\lambda$ by shifting to the neighboring pole on the
  other side if that gives a smaller value of $K_\nu$,
\item[(R2)]
if shifting to another neighboring pole is not possible,
we can invert $A-\sigma I$, where the
shift $\sigma $ is chosen near but not equal to $\lambda$, and not equal to the
neighboring poles. This results in a DPR1 matrix whose largest absolute
eigenvalue is computed accurately.
If no floating-point numbers $\sigma$ lie between $\lambda$ and the
neighboring poles, $\sigma$ and the corresponding DPR1 matrix must be computed
in double the working precision.
\end{enumerate}

We need to address one more special situation. If $\lambda$ is much closer to
zero than to the neighboring pole or poles\footnote{There can be at most one such
  eigenvalue.},
$|\lambda|\ll \min\{ |\lambda-d_{k}|,|\lambda-d_{k-1}|\}$,
then the formula (\ref{8}) may involve large cancellation, and
$\lambda$ may be inaccurate in spite of the accurately computed $\nu$ and
$v$. In this case, $\lambda$ can be computed accurately as $\lambda=1/\nu$,
where $\nu$ is the largest absolute eigenvalue of $A^{-1}$. If all poles are
non-zero, the inverse of $A$ is again an unreduced DPR1 matrix of the form
\begin{equation}\label{eq:gamma}
A^{-1}=%
D^{-1}
+\gamma D^{-1}z z^{T} D^{-1}, \quad \gamma =-\frac{\rho}{1+\rho z^{T}D^{-1}z}.
\end{equation}%
If the denominator in $\gamma$ is computed as zero, the matrix $A$ is
numerically singular and we can set $\lambda=0$.

The described procedure is implemented in algorithm $dpr1eig$.

\begin{algorithm}
\caption{}
$[\lambda,v]=\mathbf{dpr1eig}\left(
D,z,\rho ,k\right)$

\% Computes the $k$-th eigenpair
of an ordered irreducible DPR1 matrix

\% $A=\diag\left( D\right) +\rho z
z^{\prime }, \ \rho>0$

\% Find the shift $\sigma=d_i$ such that $d_i$ is the pole nearest to
  $\lambda$

\% Exterior eigenvalue $k=1$:

\textbf{if }$k==1$

\ \ $\sigma=d_{1}$

\textbf{else }

\ \ \% Interior eigenvalues $k\in\{2,\ldots,n\}$:

\ \ $\bar D=D-d_{k}$

\ \ $\tau=\bar D_{k-1}/2$

\ \ $F=1+\rho\sum (z.*z./(\bar D-\tau))$

\ \ \textbf{if }$F>0$

\ \ \ \ $\sigma=d_{k}$

\ \ \textbf{else}

\ \ \ \ $\sigma=d_{k-1}$

\ \ \textbf{end}

\textbf{end}

compute the arrowhead matrix 
$A_i^{-1}\equiv (A-\sigma I)^{-1}$
according to (\ref{invAi}) and (\ref{b_ob})

compute $\kappa_\nu$ from (\ref{knu2})

\textbf{if} $\kappa_\nu\gg O(n)$

\ \ recompute $b$ from (\ref{b_ob}) by using double the working precision
(c.f. section \ref{double})

\textbf{end}

\textbf{if} $\sigma=d_{k-1}$

\ \ compute the leftmost eigenvalue $\nu$ of $A_{i}^{-1}$ by bisection (c.f. section
\ref{sec:fast})

\textbf{else}

\ \  compute the rightmost eigenvalue $\nu$ of $A_{i}^{-1}$ by bisection

\textbf{end}

compute $v$ by (\ref{eigenvecAi}), where $\mu=1/\nu$

compute $\lambda=\mu+\sigma$

compute $K_{\nu}$ from $(\ref{K0})$

\textbf{if} $K_{\nu} \gg 1$

\ \ apply one of the remedies (R1) or (R2)

\textbf{end}

\textbf{if} $|\lambda|\ll \min\{ |\lambda-d_{k}|,|\lambda-d_{k-1}|\}$

\ \ recompute $\lambda$ from $A^{-1}$

\textbf{end}
\label{alg1}
\end{algorithm}

The algorithm $dpr1eig$ extends naturally to the Hermitian case (c.f. \cite[\S 6.1]{JSB13}).

\subsection{Accuracy of the algorithm}

Let $(\widetilde{\lambda},\widetilde{v})$ denote the eigenpair computed by
Algorithm \ref{alg1} in the standard floating-point arithmetic.
Let $\widetilde{\nu}$ denote the computed eigenvalue of $A_i^{-1}$. If $\widetilde{\nu}$ is
the absolutely largest eigenvalue of $A_i^{-1}$ and if it is computed by
bisection, then the error bound from
\cite[\S 3.1]{OlSt90} immediately implies that\footnote{Notice that a similar error bound holds for all eigenvalues which are of
the same order of magnitude as $\nu$.}
\begin{equation}  \label{lambda_max}
\widetilde{\nu}=\nu(1+\kappa_{bis}\varepsilon_{M}),\quad \kappa_{bis}\leq 1.06n\left( \sqrt{n}+1\right).
\end{equation}

The computed eigenpair satisfies
\begin{align*}
\widetilde{\lambda }&=fl(\lambda) =\lambda ( 1+\kappa
_{\lambda }\varepsilon _{M}) , \\
\widetilde{v_i }&=fl( v_i ) =v_i ( 1+\kappa _{v_i}\varepsilon _{M}) , \quad i=1,\ldots,n,
\end{align*}
where $$|\kappa_\lambda|,|\kappa_{v_i}|\leq O(\kappa_{\nu}+\kappa_{bis}),$$
and $\kappa_\nu$ is defined by (\ref{knu2}).

If $1\ll K_b\leq O(1/\varepsilon_M)$, then, after evaluating $b$ with double the
working precision, $\kappa_\nu$ is given by (\ref{knu2}) with $K_b$ replaced
by $K_b\varepsilon_M$.\footnote{If $K_b\geq O(1/\varepsilon_M)$, that is, if
  $K_b=1/\varepsilon_E$ for some $\varepsilon_E<\varepsilon_M$, then $b$ needs
  to be computed with extended precision $\varepsilon_E$.
Usage of higher precision in conjunction with the
eigenvalue computation for DPR1 matrices is analyzed in \cite{Bar93},
but there the higher precision
computation is potentially needed in the iterative part. This is less
convenient than our approach where the higher precision computation is used
only to compute one element.
}
With our approach componentwise high relative accuracy of the computed normalized eigenvectors
implies, in turn, their numerical orthogonality.

The proofs of the above error bounds are similar to the error analysis in
\cite{JSB13}.

\subsection{Fast secular equation solvers}\label{sec:fast}
Instead of using bisection to compute zeros of secular equation (\ref{secAi})
in Algorithm \ref{alg1}, we can use some fast zero finder with quadratic or even cubic
convergence like those from \cite{Mel95,BG92,Li94}. Such zero finders  compute
zeros to machine accuracy using a small number of direct evaluations of the
Pick function and its its derivative, where $O(\log(\log(1/\varepsilon)))$
iterations are needed to obtain an $\varepsilon$-accuracy \cite{LB02}.

In particular, we tested the implementation of the cubically convergent zero
finder by Borges and Gragg from \cite[\S 3.3]{BG92}, with the stopping criterion defined by
\cite[p.\ 15]{BG92}. From \cite[(21)]{BG92}, it follows that the accuracy of
the computed solution  satisfies a similar backward error bound as
(\ref{lambda_max}). This was indeed, true in all our tests. The number of
iterations never exceeded 7.

Similarly, for the solution of the secular equation (\ref{Pick}), which may be
needed in the last two ``if'' statements in Algorithm \ref{alg1}, one can
use the fast secular equation solver by Li \cite{Li94}. This solver
is implemented in the LAPACK routine \verb!DLAED4!. The accuracy of
the computed solution  satisfied a similar backward error bound as
(\ref{lambda_max}) and the number of iterations behaved as predicted.

Although the operation count of both fast zero finders is approximately half of the
operations needed for bisection, we observed no speed-up in Matlab
implementation. 

\subsection{Implementation
of the double the working precision}
\label{double}

We tried three different implementations of the double the working precision:
\begin{itemize}
\item
by converting all quantities in the formulas (\ref{b_ob}) or (\ref{eq:gamma})
to variable precision by Matlab \cite{Mat} command
\verb!sym! with parameter \verb!'f'!, and then performing the computations;
\item by evaluating all parts of the formulas (\ref{b_ob}) or
  (\ref{eq:gamma}) using extended precision routines \verb!add2!,
  \verb!sub2!, \verb!mul2!, and \verb!div2! from \cite{Dek71}; and
\item by converting all quantities in the formulas (\ref{b_ob}) or
  (\ref{eq:gamma}) from standard 64 bit double precision
numbers, declared by \verb!REAL(8)!, to 128 quadruple precision numbers,
declared by \verb!REAL(16)!, in Intel FORTRAN compiler \verb|ifort| \cite{Int}, and then
performing the computations.
\end{itemize}

Having to invoke higher precision clearly slows the
computation down.
In Matlab, when using variable precision \verb!sym! command,
the computation may be slowed down by a factor of
three hundred or more for each eigenvalue that requires formulas (\ref{b_ob})
or (\ref{eq:gamma}) to be evaluated in  higher precision. This makes use of
\verb!sym! prohibitive for higher dimensions.
Extended precision routines by Dekker
\cite{Dek71} require on average ten floating-point operations.
The fastest implementation is the one in \verb!ifort! which is only about
three times slower. Thus, the algorithm benefits from a good
implementation of higher precision.

\section{Numerical Examples}
\label{sec:ex}

We have used the following implementations of Algorithm 1:
\begin{itemize}
\item $dpr1eig(M)$ -  Matlab implementation, with double the working precision implemented
  using extended precision routines from \cite{Dek71}.
\item $dpr1eig(J)$ - Julia \cite{Julia} implementation, with double the working precision implemented
  using Julia package \verb!DoubleDouble.jl! \cite{JuliaPack} -- this
  implementation is publicly available in the Julia package
  \verb!Arrowhead.jl! \cite{JuliaPack}, and is our preferred implementation.  
\end{itemize}

We compared Algorithm 1 with the following routines:
\begin{itemize}
\item $eig$ - Matlab's standard eigenvalues routine.
\item $dlaed9$ - LAPACK routine \verb!DLAED9! compiled with \verb!ifort!
  Fortran compiler.
\item $Math$ - Mathematica \cite{Wol} eigenvalue routine with 100 digits of precision (properly
rounded to 16 decimal digits).
\end{itemize}


We illustrate our algorithm with four numerically demanding examples.
Examples 1 and 2 illustrate Algorithm 1, Example 3 illustrates the use of
double precision arithmetic, Example 4 illustrates an application to higher
dimension, and Example 5 shows the effect of using double the working precision on
overall timing. Since $dpr1eig(M)$ and $dpr1eig(J)$ give numerical identical
results, we denote these results by $dpr1eig$.

\begin{example}
\label{example1} In this example  quantities  $K_{b}$ from
(\ref{knu2}) are approximately $1$ for all eigenvalues, so we
guarantee that all eigenvalues and all components of their
corresponding eigenvectors are computed with high relative
accuracy by Algorithm \ref{alg1}, using only standard machine precision.
Let $A=D+ z z^{T}$, where
\begin{align*}
D&=\diag\, (10^{10},5,4\cdot 10^{-3},0,-4\cdot
10^{-3},-5), \\
z&=
\begin{bmatrix}
10^{10} & 1 & 1 & 10^{-7} & 1 & 1%
\end{bmatrix}^T.
\end{align*}
The computed eigenvalues are:
\footnote{If, in the last
  column, the last digits computed by {\em dpr1eig} and Mathematica,
  respectively,  differ, they are displayed in parentheses.}
{\small
\begin{equation*}
\begin{array}{lll}
\lambda ^{(eig)} & \lambda ^{(dlaed9)} & \lambda ^{(dpr1eig,Math)} \\
1.000000000100000\cdot 10^{20} & 1.000000000100000\cdot 10^{20} & 1.000000000100000\cdot 10^{20} \\
5.000000000099998 &             5.000000000100000 &             5.000000000100000 \\
4.000000099999499\cdot 10^{-3} & 4.000000100000001\cdot 10^{-3} & 4.000000100000001\cdot 10^{-3}  \\
1.665334536937735\cdot 10^{-16} &1.000000023272195\cdot 10^{-24} & 9.99999999899999(7,9) \cdot 10^{-25} \\
0 & -3.999999900000001\cdot 10^{-3} & -3.999999900000001\cdot 10^{-3} \\
 -25.00000000150000 &                -4.999999999900000 &       -4.999999999900000%
\end{array}%
\end{equation*}%
}We see that all eigenvalues
computed by $dpr1eig$ (including the tiniest ones), are exact to the machine
precision. The eigenvalues computed by $dlaed9$ are all accurate,
except $\lambda_4$. The eigenvalues computed by \verb|eig| are accurate
according to the standard perturbation theory, but they have almost no
relative accuracy\footnote{The displayed eigenvalues are the ones obtained by
the Matlab command \texttt{[V,Lambda]=eig(A)}. The command \texttt{Lambda=eig(A)}
produces slightly different eigenvalues. The reason is that Matlab uses LAPACK
routine \texttt{dsyev.f}, which, in turn, uses different algorithms depending whether
eigenvectors are required or not.}.
Due to the the accuracy of
the computed eigenvalues, the eigenvectors computed by $dpr1eig$ are
componentwise accurate up to machine precision, and therefore, orthogonal up
to machine precision. The eigenvectors computed by $dlaed9$ are also
componentwise accurate, except for
$v_4$:
{\small
\[
\begin{array}{ll}
v_4^{(dlaed9)} & v_4^{(dpr1eig,Math)}\\
  1.000000011586098\cdot 10^{-17} & 9.99999999899999(6,9)\cdot 10^{-18} \\
  2.000000023172195\cdot 10^{-18} & 1.999999999800000\cdot 10^{-18}\\
 2.500000028965244\cdot 10^{-15} & 2.499999999749999\cdot 10^{-15}\\
-1.000000000000000 & -1.000000000000000 \\
 -2.500000028965244\cdot 10^{-15} & -2.499999999749999\cdot 10^{-15} \\
  -2.000000023172195\cdot 10^{-18}& -1.999999999800000\cdot 10^{-18} %
\end{array}%
\]%
}
\end{example}

\begin{example}
\label{example2} In this example, despite very close diagonal elements, we
again guarantee that all eigenvalues and all components of their
corresponding eigenvectors are computed with high relative accuracy. Let $A=D+zz^{T}$,
where
\begin{align*}
D&=\diag\,(1+40\varepsilon ,1+30\varepsilon
,1+20\varepsilon ,1+10\varepsilon ), \\
z&=
\begin{bmatrix}
1 & 2 & 2 & 1
\end{bmatrix}.
\end{align*}
and $\varepsilon=2^{-52}=2\varepsilon _{M}$. For
this matrix, the quantities $K_{b}$ are again of order one for all
eigenvalues, so Algorithm \ref{alg1} uses only standard working precision.
The computed eigenvalues are:
{\small
\[
\begin{array}{lll}
\lambda^{(eig)} & \lambda^{(dlaed9)} & \lambda^{(dpr1eig)} \\
11+32\varepsilon  & 11+48\varepsilon  & 11+32\varepsilon  \\
1+38\varepsilon  & 1+41\varepsilon  & 1+39\varepsilon  \\
1+31\varepsilon  & 1+27\varepsilon  & 1+25\varepsilon  \\
1+8\varepsilon  & 1+9\varepsilon  & 1+11\varepsilon
\end{array}
\]
}Notice that all computed eigenvalues are accurate according to standard
perturbation theory. However, only the eigenvalues computed by $dpr1eig$
satisfy the interlacing property.
The eigenvalues computed by $Math$, properly rounded to 32 decimal digits are:
{\small
\[
\begin{array}{l}
\lambda ^{(Math)} \\
 11.000000000000005551115123125783\\
 1.0000000000000085712482686374087\\
 1.0000000000000055511151231257826\\
 1.0000000000000025309819776141565\\
\end{array}
\]
}If Algorithm \ref{alg1} is modified to
return $\sigma$ and $\mu$ (both in standard precision), then for the
eigenvalues $\lambda_2$, $\lambda_3$ and $\lambda_4$ the corresponding pairs
$(\sigma,\mu)$ give representations of those eigenvalues to 32 decimal digits.
In our case, the exact values $\sigma+\mu$ properly rounded to 32 decimal digits are
equal to the corresponding eigenvalues computed by Mathematica displayed
above.

The eigenvectors $v_2$, $v_3$ and $v_4$ computed by $eig$ span an invariant
subspace of $\lambda_2$, $\lambda_3$ and $\lambda_4$, but their components are not accurate.
Due to the accuracy of the computed eigenvalues, the eigenvectors computed by
$dpr1eig$ are componentwise accurate up to the machine precision (they coincide with
the eigenvectors computed by $Math$, and are therefore
orthogonal. Interestingly, in this example the eigenvectors computed by
$dlaed9$ are also componentwise accurate, but there is no underlying
theory for such high accuracy.

\end{example}

\begin{example}
\label{example3} In this example (see \cite{GE94}) we can
guarantee that all eigenvalues and eigenvectors will be computed with
componentwise high relative accuracy only if $b$ from (\ref{b_ob})
is for $k\in\{2,3,4\}$ computed in double of the working precision.
Let $A=D+zz^{T}$,
where
\begin{align*}
D&=\diag\,(10/3,2+\beta ,2-\beta ,1),\\
z&=\begin{bmatrix}
2 & \beta & \beta & 2
\end{bmatrix}, \quad \beta =10^{-7}.
\end{align*}

For $k\in\{2,3,4\}$ the quantities $\kappa_{\nu}$ from (\ref{knu2}) are of
order $O(10^7)$, so the element $b$ in each of the
matrices needs to be computed in double of the working precision.
For example, for $k=2$, the element $b=\left[
A_{2}^{-1}\right] _{22}$ computed by Algorithm \ref{alg1} in standard precision is
equal to $b= 5.749999751891721\cdot 10^{7}$, while Matlab routine
\verb|inv| gives $b=5.749999746046776\cdot 10^{7}$.
Computing $b$ in double of the working precision in Algorithm \ref{alg1} gives the
correct value $b=5.749999754927588\cdot 10^{7}$.

The eigenvalues computed by $eig$, $dlaed9$, $dpr1eig$ and $Math$, respectively, are all highly
relatively accurate -- they differ in the last or last two digits.
However, the eigenvectors $v_2$, $v_3$ and $v_4$ computed by $dpr1eig$ (with
double precision computation of $b$'s), are componentwise accurate to machine
precision and therefore orthogonal. 
The eigenvectors computed by $eig$ and $dlaed9$ are, of course, orthogonal,
but are not componentwise accurate. For
example,
\[
\begin{array}{lll}
 v_2^{(eig)} & v_2^{(dlaed9)} &v_2^{(dpr1eig,Math)} \\
     2.088932176072975\cdot 10^{-1}&  2.088932143122528\cdot 10^{-1}  &
     2.088932138163857\cdot 10^{-1}\\
    -9.351941376557037\cdot 10^{-1}&  -9.351941395201120\cdot 10^{-1} &
    -9.351941398441738\cdot 10^{-1}\\
    -6.480586028358029\cdot 10^{-2}&   -6.480586288204153\cdot 10^{-2} &
    -6.480586264549802\cdot 10^{-2}\\
    -2.785242341430628\cdot 10^{-1}& -2.785242297496694\cdot 10^{-1} &
    -2.785242290885133\cdot 10^{-1}
\end{array}
\]
\end{example}

\begin{example}
\label{example4}
In this example we extend Example \ref{example3} to higher dimension, as
in TEST 3 from \cite[\S 6]{GE94}.
Here $A=D+zz^{T}\in\mathbb{R}^{202\times 202}$, where
\begin{align*}
D&=\diag\,(1,2+\beta ,2-\beta,2+2\beta ,2-2\beta,\ldots,2+100\beta,2-100\beta, 10/3),\\
z&=\begin{bmatrix}
2 & \beta & \beta & \ldots & \beta & 2
\end{bmatrix}, \quad \beta \in\{10^{-3},10^{-8},10^{-15}\}.
\end{align*}
For each $\beta$, we solved the eigenvalue problem with Algorithm \ref{alg1} without using
double the working precision ($dpr1eig\_nd$), $dpr1eig$,
and $dlaed9$.
For $\beta=10^{-3}$, Algorithm \ref{alg1} used double the working precision for
computing 25 eigenvalues, and for $\beta=10^{-8}$ and
$\beta=10^{-15}$ double the working precision was needed for all but the largest
eigenvalue. As in \cite[\S 6]{GE94}, for each algorithm we computed orthogonality
and residual measures,
$$
\mathcal{O}=\max_{1\leq i\leq n}\frac{\| V^T
  v_i-e_i\|_2}{n\varepsilon_M},\qquad
\mathcal{R}=\max_{1\leq i\leq n}\frac{\| A
  v_i-\lambda_iv_i\|_2}{n\varepsilon_M\|A\|_2},
$$
respectively. Here $V=\begin{bmatrix} v_1&v_2 &\cdots &v_n
\end{bmatrix}$ is the computed matrix of eigenvectors, and $e_i$ is the $i$-th
column of the identity matrix.


Since we proved the componentwise accuracy of
eigenvectors computed by $dpr1eig$, we take those as the ones of reference.
Table \ref{tab:1} displays orthogonality measures, residual measures, relative
errors in the computed eigenvalues and componentwise relative errors in the
computed eigenvectors, superscripted by the name of the respective algorithm.
\begin{table}[hbtp]
  \centering
  \begin{tabular}{cccc}
    $\beta$ & $10^{-3}$ & $10^{-8}$ & $10^{-15}$ \\  \hline\hline
    $\mathcal{O}^{(dpr1eig\_nd)}$ & 1.47  & $5.8\cdot 10^4$  &  $2.1\cdot 10^{11}$   \\
    $\mathcal{O}^{(dpr1eig)}$ &  0.059    & 0.039  & 0.045    \\
    $\mathcal{O}^{(dlaed9)}$ &  0.049    &  0.064 & 0.045   \\ \hline
    $\mathcal{R}^{(dpr1eig\_nd)}$ & 0.0086  & 0.033  & 0.0043    \\
    $\mathcal{R}^{(dpr1eig)}$ &  0.0086    & 0.039  & 0.0043    \\
    $\mathcal{R}^{(dlaed9)}$ &  0.029    &  0.03  &  0.013 \\ \hline
 $\max\limits_{1\leq i\leq n}\frac{| \lambda^{(dpr1eig\_nd)}_i-\lambda^{(dpr1eig)}_i|}
 {|\lambda^{(dpr1eig)}_i|}$  &  $2.2\cdot 10^{-16}$      &  0  & $2.2\cdot 10^{-16}$   \\
 $\max\limits_{1\leq i\leq n}\frac{| \lambda^{(dlaed9)}_i-\lambda^{(dpr1eig)}_i|}
 {|\lambda^{(dpr1eig)}_i|}$  &  $1.5\cdot 10^{-15}$      &  $2.2\cdot 10^{-16}$     & 0  \\ \hline
 $\max\limits_{1\leq i,j\leq n}\frac{| [v^{(dpr1eig\_nd)}_i]_j- [v^{(dpr1eig)}_i]_j|}
 {|[v^{(dpr1eig)}_i]_j|}$  & $2.7\cdot 10^{-13}$       & $2.8\cdot 10^{-8}$   &  0.518  \\
$\max\limits_{1\leq i,j\leq n}\frac{| [v^{(dlaed9)}_i]_j- [v^{(dpr1eig)}_i]_j|}
 {|[v^{(dpr1eig)}_i]_j|}$  & $2.2\cdot 10^{-12}$       & $1.9\cdot 10^{-8}$   & 0.043   \\ \hline
  \end{tabular}
  \caption{Orthogonality measures, residue measures, relative errors in computed eigenvalues, and
    componentwise relative errors in computed eigenvectors.}
\label{tab:1}
\end{table}
From table \ref{tab:1}, we see that all algorithms behave exactly as
predicted by the theoretical analysis. All algorithms compute all eigenvalues
to high relative accuracy  because it is the same as normwise accuracy for
this case.
$dpr1eig\_nd$ loses orthogonality as predicted by the respective
condition numbers. The number of correct digits in the computed eigenvectors
is approximately the same for $dpr1eig\_nd$ and $dlaed9$, but
there is no proof of such componentwise accuracy of the eigenvectors computed
by $dlaed9$.
As a consequence of their componentwise accuracy, the eigenvectors computed by
$dpr1eig$ are fully orthogonal. 
\end{example}

\begin{example}\label{example5}
To illustrate the effect of using double the working precision, in Table
\ref{tab:2} we give
timings for the matrix $A\in \mathbb{R}^{202\times 202}$ of the same form as
in Example 4.

\begin{table}[hbtp]
  \centering
  \begin{tabular}{c|ccc}
 & $\beta=10^{-3}$  & $\beta=10^{-8}$  & $\beta=10^{-15}$ \\ \hline
$dpr1eig(M)$ &  11 & 17 & 17   \\
$dpr1eig(J)$ &  1.2 & 2.2 & 2.2 \\
$dlaed9$     &  0.13 & 0.13 & 0.13
  \end{tabular}
  \caption{Running time (in seconds) for the computation of eigenvalues and
    eigenvectors of DPR1 matrix $A$ of order $n=2002$.}
\label{tab:2}
\end{table}

We see that the Julia version of Algorithm 1 is almost 10 times faster than
the Matlab version, which makes Julia version  an implementation of preference. 
As in Example 4, for $\beta=10^{-3}$, $dpr1eig$ used double the working
precision to compute respective $b$ when
computing 25 eigenvalues, and for $\beta=10^{-8}$ and
$\beta=10^{-15}$ double the working precision was needed for all but the largest
eigenvalue. We see that the overhead of using double the working precision 
is approximately 55\% in both, Julia and Matlab. 
  
\end{example}

\section*{Acknowledgment}
We would like to thank Ren Cang Li for providing Matlab implementation of the
LAPACK routine \verb!DLAED4! and its dependencies.

\end{document}